\documentclass[11pt]{article}
\usepackage{enumerate}
\usepackage{amssymb,a4wide,latexsym,makeidx,epsfig,fleqn}
\usepackage{amsthm}
\usepackage{amsmath}
\usepackage{enumerate}
\newtheorem{theorem}{Theorem}[section]

\newtheorem{definition}[theorem]{Definition}
\newtheorem{lemma}[theorem]{Lemma}

\begin{document}
\textwidth 150mm \textheight 225mm
\title{Skew-rank of an oriented graph in terms of the rank  and dimension of cycle space of its underlying graph \footnote{This work is supported by the National Natural Science Foundation of China (No. 11171273).}}
\author{{Yong Lu, Ligong Wang\footnote{Corresponding author.} ~and Qiannan Zhou}\\
{\small Department of Applied Mathematics, School of Science, Northwestern
Polytechnical University,}\\ {\small  Xi'an, Shaanxi 710072,
People's Republic
of China.}\\
{\small E-mail: luyong.gougou@163.com, lgwangmath@163.com, qnzhoumath@163.com.}}

\date{}
\maketitle
\begin{center}
\begin{minipage}{120mm}
\vskip 0.3cm
\begin{center}
{\small {\bf Abstract}}
\end{center}
{\small Let $G^{\sigma}$ be an oriented graph and $S(G^{\sigma})$ be its skew-adjacency matrix, where $G$ is called the underlying graph of  $G^{\sigma}$. The skew-rank of $G^{\sigma}$, denoted by $sr(G^{\sigma})$, is the rank of $S(G^{\sigma})$. Denote by $d(G)=|E(G)|-|V(G)|+\theta(G)$ the dimension of cycle spaces of $G$, where $|E(G)|$, $|V(G)|$ and $\theta(G)$ are the edge number, vertex number and the number of connected components of $G$, respectively.
Recently, Wong, Ma and Tian  [European J. Combin. 54 (2016) 76--86] proved that $sr(G^{\sigma})\leq r(G)+2d(G)$ for an oriented graph $G^{\sigma}$, where $r(G)$ is the rank of the adjacency matrix of $G$, and characterized the graphs whose skew-rank attain the upper bound. However, the problem of the lower bound of $sr(G^{\sigma})$ of an oriented graph $G^{\sigma}$ in terms of $r(G)$ and $d(G)$ of its underlying graph $G$ is left open till now. In this paper, we prove that  $sr(G^{\sigma})\geq r(G)-2d(G)$ for an oriented graph $G^{\sigma}$ and characterize the graphs whose skew-rank attain the lower bound.

\vskip 0.1in \noindent {\bf Key Words}: \ Skew-rank,  Rank of graphs, Dimension of cycle space. \vskip
0.1in \noindent {\bf AMS Subject Classification (2010)}: \ 05C50. }
\end{minipage}
\end{center}

\section{Introduction}

In this paper, we only consider simple graphs without multiple edges and loops.
Let $G$ be a  simple graph with vertex set $V(G)=\{v_{1},v_{2},\ldots,v_{n}\}$ and edge set $E(G)$. The \emph{adjacency matrix} of  $G$  of order $n$ is defined as the $n\times n$ symmetric square matrix $A=A(G)=(a_{ij})$, where  $a_{ij}=1$  if $v_{i}v_{j}\in E(G)$, otherwise $a_{ij}=0$. The \emph{rank} $r(G)$ of $G$ is defined to be the rank of $A(G)$, and the \emph{nullity} $\eta(G)$ of $G$ is defined to be the multiplicity of 0 as an eigenvalue of $A(G)$. Obviously, $|V(G)|=r(G)+\eta(G)$.
We use Bondy and Murty \cite{BM} for terminologies and notations not defined here.

An oriented graph $G^{\sigma}$ is a digraph which assigns each edge of $G$ with a direction $\sigma$, where $G$ is called the \emph{underlying graph} of $G^{\sigma}$. The \emph{skew-adjacency matrix} associated to $G^{\sigma}$ is the $n\times n$ matrix $S(G^{\sigma})=(s_{ij})$, where $s_{ij}=-s_{ji}=1$  if $(v_{i},v_{j})$ is an arc of $G^{\sigma}$, otherwise $s_{ij}=s_{ji}=0$. The \emph{skew-rank} $sr(G^{\sigma})$ of an oriented graph $G^{\sigma}$ is defined as the rank of the skew-adjacency matrix $S(G^{\sigma})$. Since $S(G^{\sigma})$ is skew-symmetric, every eigenvalue of $S(G^{\sigma})$ is a pure imaginary number or 0, and the skew-rank of an oriented graph is even.

Let $C_{n}^{\sigma}=v_{1}v_{2}\cdots v_{n}v_{1}$ be an even oriented cycle. Denote by sgn($C_{n}^{\sigma}$) the \emph{sign}  of $C_{n}^{\sigma}$, which is defined as the sign of $\prod_{i=1}^{n}s_{v_{i}v_{i+1}}$ with $v_{n+1}=v_{1}$. An even oriented cycle $C_{n}^{\sigma}$ is called \emph{evenly-oriented}  (resp., \emph{oddly-oriented}) if its sign is positive (resp., negative). $G^{\sigma}$ is called \emph{evenly-oriented} if every even cycle in $G^{\sigma}$ is evenly-oriented.

Sometimes we use the notation $G^{\sigma}-H^{\sigma}$ instead of $G^{\sigma}-V(H^{\sigma})$ if $H^{\sigma}$ is an \emph{induced subgraph} of $G^{\sigma}$, where $G^{\sigma}-H^{\sigma}$ is the subgraph obtained from $G^{\sigma}$ by deleting all vertices of $H^{\sigma}$ and all incident edges. For an induced subgraph $H^{\sigma}$ and a vertex $x$ outside $H^{\sigma}$, the induced subgraph of $G^{\sigma}$ with vertex set $V(H^{\sigma})\cup\{x\}$ is simply written as $H^{\sigma}+x$.  For a vertex $v\in V(G^{\sigma})$, let $G^{\sigma}-v$ denote the oriented graph obtained from $G^{\sigma}$ by removing the vertex $v$ and all edges incident with $v$. A vertex $x\in V(G)$ is called a \emph{cut-point} of a connected graph $G$ if the resultant graph $G-x$ has at least two components. A vertex of $G^{\sigma}$ is called a \emph{pendant} vertex if its degree is 1 in $G$, and is called a \emph{quasi-pendant} vertex if it is adjacent to a pendant vertex.  An induced subgraph $H$ (resp., $H^{\sigma}$) of a graph $G$ (resp., $G^{\sigma}$) is called a \emph{pendant cycle} (resp.,  \emph{pendant oriented cycle}) of $G$ (resp., $G^{\sigma}$) if $H$ is a cycle such that $H$ has a unique vertex of degree 3 in $G$. Denote by $d(G)=|E(G)|-|V(G)|+\theta(G)$ the \emph{dimension of cycle spaces} of $G$, where $|E(G)|$, $|V(G)|$ and $\theta(G)$ are the edge number, vertex number and the number of connected components of $G$, respectively. Obviously, when $G$ is connected, then $G$ is a tree if $d(G)=0$, and $G$ is a unicyclic graph if $d(G)=1$.
A \emph{matching} in a graph $G$ is a set of pairwise nonadjacent edges. A  \emph{maximum matching} is one that contains as many edges of $G$ as possible. The \emph{matching number} of $G$, denoted by $m(G)$, is the size of a maximum matching in $G$.
Denote by $P_{n},~C_{n}$ a path and a cycle  of order $n$, respectively. A graph is called \emph{empty} if it has some vertex and no edges.

In 1957, Collatz and Sinogowitz \cite{CS} first posed the problem of characterizing all graphs $G$ with $\eta(G)>0$. This problem is of great interest in both chemistry and mathematics. For a bipartite graph $G$ which corresponds to an alternant hydrocarbon in chemistry, if $\eta(G)>0$, it is indicated that the corresponding molecule is unstable. The nullity of a graph is also meaningful in mathematics since it is related to the singularity of adjacency matrix.

Till now, many scholars investigated the nullity of graphs, they focused on special graph classes, such as trees, unicyclic graphs, bicyclic graphs, bipartite graphs and so on. There are also some papers focused on the study of the connection between the nullity (or rank) of graphs $G$ in terms of certain structural parameters, such as matching number, dimension of cycle spaces and so on.
Recently, Wang and Wong \cite{WW} obtained the  bounds for the matching number of $G$ in terms of the $r(G)$ and $d(G)$, that is: $$\left\lceil\frac{r(G)-d(G)}{2}\right\rceil\leq m(G)\leq\left\lfloor\frac{r(G)+2d(G)}{2}\right\rfloor.$$
The bounds for the matching number can be rewritten in an equivalent form as bounds for the nullity of $G$, that is:$$|V(G)|-2m(G)-d(G)\leq\eta(G)\leq|V(G)|-2m(G)+2d(G).$$
In 2015, Song, Song and Tam \cite{SST} characterized the graph $G$ that satisfy the equality $\eta(G)=|V(G)|-2m(G)+2d(G)$. The lower bound $|V(G)|-2m(G)-d(G)$ of $\eta(G)$  was characterized by Wang \cite{LWANG} and independently by Rula, Chang and Zheng in \cite{RCZ}.

In 2016, Ma, Wong and Tian \cite{MWTDAM}  proved that $$\eta(G)\leq2d(G)+p(G),$$ where $p(G)$ is the number of pendant vertices of $G$, they also proved that the equality is attained if and only if every component of $G$ is a cycle with size a multiple of 4.

Recently, the skew-rank of  skew-adjacency matrix of an oriented graph has received a lot of attentions.
Li and Yu \cite{LXYG} studied the skew-rank of oriented graphs and characterized oriented unicyclic graphs  attaining the minimum value of the skew-rank among oriented unicycle graphs of order $n$ with girth $k$. Qu and Yu \cite{QHYG} characterized the bicyclic oriented graphs with skew-rank 2 or 4. Lu, Wang and  Zhou \cite{LWZ} characterized the bicyclic oriented graphs with skew-rank 6. Qu, Yu and Feng \cite{QHYGFL} obtained more results about the minimum skew-rank of graphs. They also characterized the unicyclic graphs with skew-rank 4 or 6, respectively.

In \cite{CT}, Chen and Tian proved that $sr(G^{\sigma})\geq\sum_{i=1}^{k}q_{i}-2k$ if $G^{\sigma}$ is a connected oriented graph with $k$ pairwise edge-disjoint cycles of orders $q_{1},q_{2},\ldots,q_{k}$.  In \cite{MWTLAA}, Ma, Wong and Tian characterized the bounds of skew-rank of an oriented connected graph $G^{\sigma}$ in terms of matching number, that is: $$2m(G)-2\beta(G)\leq sr(G^{\sigma})\leq2m(G),$$ where $\beta(G)=|E(G)|-|V(G)|+1$. The oriented graphs satisfying $sr(G^{\sigma})=2m(G)-2\beta(G)$ are characterized definitely.

In 2016, Wong, Ma and Tian \cite{WMT} proved that $$sr(G^{\sigma})\leq r(G)+2d(G)$$ for an oriented graph $G^{\sigma}$. They characterized the oriented graphs $G^{\sigma}$ whose skew-rank can attains the upper bound.

A natural problem is : How about the lower bound of the skew-rank of an oriented graph $G^{\sigma}$ in terms of the rank  and the dimension of cycle spaces $d(G)$ of its underlying graph $G$? In this paper, we will prove that $$sr(G^{\sigma})\geq r(G)-2d(G)$$ for an oriented graph $G^{\sigma}$ and characterize the oriented graphs $G^{\sigma}$ whose skew-rank can attains the lower bound. Our main results are Theorems \ref{th:1.1} and \ref{th:1.3}.

\noindent\begin{theorem}\label{th:1.1}
Let $G^{\sigma}$ be a finite oriented graph without loops and multiple arcs. Then $$sr(G^{\sigma})\geq r(G)-2d(G).$$
\end{theorem}

Combining with the upper bound of the skew-rank of an oriented graph $G^{\sigma}$ in \cite{WMT} and our result in Theorem \ref{th:1.1}, we have $$r(G)-2d(G)\leq sr(G^{\sigma})\leq r(G)+2d(G).$$

\noindent\begin{definition}\label{de:1.2} (\cite{WMT})
\rm Let $G$ be a graph with at least one pendant vertex. The operation of deleting a pendant vertex and its adjacent vertex from $G$ is called \emph{$\delta$-transformation}.
\end{definition}

An oriented graph $G^{\sigma}$ will be called \emph{lower-optimal} if the skew-rank of $G^{\sigma}$ attains the lower bound $r(G)-2d(G)$.
A graph $G$ is called \emph{pairwise vertex-disjoint} if distinct cycles (if any) of $G$ have no common vertices.
\noindent\begin{theorem}\label{th:1.3}
Let $G^{\sigma}$ be a finite oriented graph without loops and multiple arcs of order $n$. Then $G^{\sigma}$ is lower-optimal if and only if the following conditions all hold.
\begin{enumerate}[(1)]
  \item Cycles (if any) of $G^{\sigma}$ are pairwise vertex-disjoint.
  \item Each cycle $C_{q}^{\sigma}$ of $G^{\sigma}$ is evenly-oriented  with order $q\equiv 2(\emph{mod}~4)$.
  \item A series of $\delta$-transformations can switch $G$ to a crucial subgraph $G_{0}$, which is the disjoint union of $d(G)$ cycles together with some isolated vertices.
\end{enumerate}
\end{theorem}

The rest of this paper is organized as follows: in Section 2, some necessary lemmas are introduced. In Section 3, we will prove Theorem \ref{th:1.1}.  In Section 4, we will give some useful lemmas and theorems, and prove  Theorem \ref{th:1.3}.

\section{Preliminaries}

In this section, we introduce some elementary lemmas and known results.
\noindent\begin{lemma}\label{le:2.1}(\cite{LXYG})
\begin{enumerate}[(a)]
  \item Let $H^{\sigma}$ be an induced subgraph of $G^{\sigma}$. Then $sr(H^{\sigma})\leq sr(G^{\sigma})$.
  \item Let $G^{\sigma}=G_{1}^{\sigma}\cup G_{2}^{\sigma}\cup\cdots\cup G_{t}^{\sigma}$, where $G_{1}^{\sigma},G_{2}^{\sigma},\ldots,G_{t}^{\sigma}$ are connected components of $G^{\sigma}$. Then $sr(G^{\sigma})=\sum_{i=1}^{t}sr(G_{i}^{\sigma})$.
  \item Let $G^{\sigma}$ be an oriented graph on $n$ vertices. Then $sr(G^{\sigma})=0$ if and only if $G^{\sigma}$ is an empty graph.
\end{enumerate}
\end{lemma}

Note that the results of Lemma \ref{le:2.1} also hold for the underlying graph $G$ of $G^{\sigma}$.

\noindent\begin{lemma}\label{le:2.2}(\cite{LXYG})
Let $T^{\sigma}$ be an oriented acyclic graph with matching number $m(T)$. Then $r(T)=sr(T^{\sigma})=2m(T)$.
\end{lemma}

\noindent\begin{lemma}\label{le:2.3}(\cite{LXYG})
Let $C_{n}^{\sigma}$ be an oriented cycle of order n. Then we have
\begin{displaymath}
sr(C_{n}^{\sigma})=\left\{
        \begin{array}{ll}
          n, &C_{n}^{\sigma}~ \rm is~oddly$-$\emph{oriented}, \\
          n-2,& C_{n}^{\sigma} ~\rm is~evenly$-$\emph{oriented}, \\
          n-1, &\rm otherwise. \\
        \end{array}
      \right.
\end{displaymath}
\end{lemma}


\noindent\begin{lemma}\label{le:2.4}(\cite{CL})
Let $G$ be a graph containing a pendant vertex, and $H$ be the induced subgraph of $G$ obtained by deleting this pendant vertex together with the vertex adjacent to it. Then $r(G)=r(H)+2$.
\end{lemma}

\noindent\begin{lemma}\label{le:2.5}(\cite{LXYG})
Let $G^{\sigma}$ be an oriented graph containing a pendant vertex, and $H^{\sigma}$ be the induced subgraph of $G^{\sigma}$ obtained by deleting this pendant vertex together with the vertex adjacent to it. Then $sr(G^{\sigma})=sr(H^{\sigma})+2$.
\end{lemma}

\noindent\begin{lemma}\label{le:2.6}(\cite{GX})
Let $x$ be a cut-point of a graph $G$ and $G_{1},G_{2},\ldots,G_{t}$ be all  components of $G-x$. If there exists a component, say $G_{1}$, such that $r(G_{1})=r(G_{1}+x)-2$, then $r(G)=r(G-x)+2$. If $r(G_{1})=r(G_{1}+x)$, then $r(G)=r(G_{1})+r(G-G_{1})$.
\end{lemma}

\noindent\begin{lemma}\label{le:2.7}(\cite{WMT})
Let $x$ be a vertex of $G^{\sigma}$. Then $sr(G^{\sigma}-x)$ is equal either to $sr(G^{\sigma})$ or to $sr(G^{\sigma})-2$.
\end{lemma}

\noindent\begin{lemma}\label{le:2.8}(\cite{BBD})
If $x$ is a vertex of a graph $G$, then $r(G)-2\leq r (G-x)\leq r(G)$.
\end{lemma}

\noindent\begin{lemma}\label{le:2.9}(\cite{CL})
Let $C_{q}$ be a cycle of order $q$. Then $r(C_{q})=q-2$ if $q\equiv0(\emph{mod}~4)$, and $r(C_{q})=q$ otherwise. Let $P_{q}$ be a path of order $q$, then  $r(P_{q})=q$ if $q$ is even, and $r(P_{q})=q-1$ if $q$ is odd.
\end{lemma}

\section{Proof for Theorem \ref{th:1.1}}
In this section, we will prove Theorem \ref{th:1.1}. First, we will introduce the following lemma that will be useful for later.
\noindent\begin{lemma}\label{le:3.1}(\cite{WMT})
Let $G$ be a  graph with a vertex $x$. Then
\begin{enumerate}[(a)]
  \item  $d(G)=d(G-x)$ if $x$ lies outside any cycle of $G$.
  \item $d(G-x)\leq d(G)-1$ if $x$ lies on a cycle of $G$.
  \item $d(G-x)\leq d(G)-2$ if $x$ is a common vertex of distinct cycles of $G$.
  \item If the cycles of $G$ are pairwise vertex-disjoint, then $d(G)$ precisely equals the number of cycles in $G$.
\end{enumerate}
\end{lemma}

From \cite{WMT}, we know that a similar result as Lemma \ref{le:3.1} holds for an oriented graph $G^{\sigma}$. Now, we will prove Theorem \ref{th:1.1}.\\

\noindent \textbf{Proof of Theorem \ref{th:1.1}.} We shall apply induction on $d(G)$ to prove $sr(G^{\sigma})\geq r(G)-2d(G)$.

\textbf{Case 1.} If $d(G)=0$, then the result follows from Lemma \ref{le:2.2}.

\textbf{Case 2.} If $d(G)\geq1$, then $G^{\sigma}$ has at least one cycle. Let $x$ be a vertex of a cycle of $G^{\sigma}$. By Lemma \ref{le:3.1},
\begin{equation}\label{1}
d(G-x)\leq d(G)-1.
\end{equation}
The induction hypothesis to $G^{\sigma}-x$ allows us to assume
\begin{equation}\label{2}
sr(G^{\sigma}-x)\geq r(G-x)-2d(G-x).
\end{equation}
By Lemmas \ref{le:2.1} and \ref{le:2.8},
\begin{equation}\label{3}
sr(G^{\sigma})\geq sr(G^{\sigma}-x),~r(G-x)\geq r(G)-2.
\end{equation}
Combining with inequalities (\ref{1})--(\ref{3}), we have
\begin{equation}\label{4}
sr(G^{\sigma})\geq sr(G^{\sigma}-x)\geq r(G-x)-2d(G-x)\geq r(G)-2-2d(G)+2=r(G)-2d(G).
\end{equation}

This completes the proof. \quad $\square$

Combining with the upper bound of the skew-rank of an oriented graph $G^{\sigma}$ in \cite{WMT} and our result in Theorem \ref{th:1.1}, we have $$r(G)-2d(G)\leq sr(G^{\sigma})\leq r(G)+2d(G).$$
Now, we will prove the Theorem \ref{th:1.3} to characterize the graphs whose skew-rank can attain the lower bound.

\section{Proof for Theorem \ref{th:1.3}}

In this section, we will give some useful lemmas and theorems, and prove the Theorem \ref{th:1.3}.

\noindent\begin{lemma}\label{le:4.1}
Let $x$ be a  vertex lying on a cycle of $G^{\sigma}$. If $G^{\sigma}$ is lower-optimal, then
\begin{enumerate}[(a)]
  \item  $sr(G^{\sigma})=sr(G^{\sigma}-x)$,~$r(G)=r(G-x)+2$,~$d(G)=d(G-x)+1$.
  \item $G^{\sigma}-x$ is lower-optimal.
  \item $x$ lies on only one cycle of $G$ and $x$ is not a quasi-pendant vertex of $G$.
\end{enumerate}
\end{lemma}

\noindent\textbf{Proof.}
From Theorem \ref{th:1.1} and $G^{\sigma}$ is lower-optimal, we have
$r(G)-2d(G)=sr(G^{\sigma})\geq r(G)-2d(G)$, which forces inequalities (\ref{1})--(\ref{3}) in the proof of Theorem \ref{th:1.1}, all turn into equalities. So, (a) and (b) of this lemma are all derived.

By Lemma \ref{le:3.1} and (a) of this lemma, we know that $x$ cannot be a common vertex of two distinct cycles in $G^{\sigma}$.

Suppose that $x$ is a quasi-pendant vertex adjacent to a pendant vertex $y$, by Lemma \ref{le:2.5}, we have $sr(G^{\sigma})=sr(G^{\sigma}-x-y)+2=sr(G^{\sigma}-x)+2$, which contradicts to (a) of this lemma.

This completes the proof. \quad $\square$

From Lemma 2.6 of \cite{LWZ} and Lemma 4.3 of \cite{WMT}, we have the following lemma.

\noindent\begin{lemma}\label{le:4.2}(\cite{LWZ})
Let $C_{q}^{\sigma}$ be a pendant oriented cycle of $G^{\sigma}$ with $x$ the unique vertex of $C_{q}$ of degree 3, and let $H^{\sigma}=G^{\sigma}-C_{q}^{\sigma}$, $K^{\sigma}=H^{\sigma}+x$. Then
\begin{displaymath}
sr(G^{\sigma})=\left\{\
        \begin{array}{ll}
          q-2+sr(K^{\sigma}),& C_{q}^{\sigma}\rm~is~evenly$-$\rm{oriented,}\\
          q+sr(H^{\sigma}), & C_{q}^{\sigma}\rm~is~oddly$-$\rm{oriented,}\\
          q-1+sr(K^{\sigma}),& \rm otherwise.
        \end{array}
      \right.
\end{displaymath}
\end{lemma}

\noindent\begin{theorem}\label{th:4.3}
Let $C_{q}^{\sigma}$ be a pendant oriented cycle of $G^{\sigma}$ with $x$ the unique vertex of $C_{q}$ of degree 3, and let $H^{\sigma}=G^{\sigma}-C_{q}^{\sigma}$, $K^{\sigma}=H^{\sigma}+x$. If $G^{\sigma}$ is lower-optimal, then
\begin{enumerate}[(a)]
  \item  $q\equiv 2(\emph{mod}~4)$ and $C_{q}^{\sigma}$ is evenly-oriented.
  \item $s(G^{\sigma})=q-2+sr(K^{\sigma})$, $sr(H^{\sigma})=sr(K^{\sigma})$, $r(G)=q+r(K)$ and $r(H)=r(K)$.
  \item Both $H^{\sigma}$ and $K^{\sigma}$ are lower-optimal.
\end{enumerate}
\end{theorem}
\noindent\textbf{Proof.} Assertion (a) of this theorem will be derived after  three claims.

\textbf{Claim 1.} $q$ is even.

Suppose that $q$ is odd, by Lemma \ref{le:4.2},
\begin{equation}\label{5}
sr(G^{\sigma})=q-1+sr(K^{\sigma}).
\end{equation}

Further, since $G^{\sigma}$ is lower-optimal,
$r(G)=sr(G^{\sigma})+2d(G)=q-1+sr(K^{\sigma})+2d(G)\geq q-1+r(K)-2d(K)+2d(G)=q-1+r(K)+2=q+1+r(K)$, where the inequality follows from Theorem \ref{th:1.1}.

Since $x$ lies on the cycle $C_{q}$, by (a) of Lemma \ref{le:4.1}, we have
\begin{equation}\label{6}
r(G)=r(G-x)+2=q-1+r(H)+2=q+1+r(H).
\end{equation}

So, $r(G)=q+1+r(H)\geq q+1+r(K)$, i.e., $r(H)\geq r(K)$.

From Lemma \ref{le:2.8}, we know that $r(H)\leq r(K)$. So,
\begin{equation}\label{7}
 r(H)=r(K).
\end{equation}

Let $A(G)$ be the adjacency matrix of $G$, where
\begin{displaymath}
 A(G)=\left(
  \begin{array}{cccccccccccccc}
            A&     \alpha&   0 \\
            \alpha^{T}&     0&   \beta\\
            0&    \beta^{T}&   B\\

  \end{array}
\right),
\end{displaymath}
where $A$ is the adjacency matrix of $C_{q}-x$, $B$ is the adjacency matrix of $H$, $\alpha^{T}$ refers to the transpose of $\alpha$.
From the process of the proof in Lemma 4.4 in \cite{WMT}, we have
 \begin{displaymath}
 r(G)=r\left(
  \begin{array}{cccccccccccccc}
            A&     0&   0 \\
            0&     a&   \beta\\
            0&    \beta^{T}&   B\\

  \end{array}
\right),
\end{displaymath}
where $a=-\alpha^{T}A^{-1}\alpha$.
So, \begin{displaymath}
 r(G)=r(A)+r\left(
  \begin{array}{cccccccccccccc}
               a&   \beta\\
               \beta^{T}&   B\\
  \end{array}
\right)\leq r(A)+r\left(
  \begin{array}{cccccccccccccc}
            a&         0\\
            0&          0 \\
  \end{array}
  \right)+r\left(
  \begin{array}{cccccccccccccc}
             0&   \beta\\
               \beta^{T}&   B\\
  \end{array}
  \right).
\end{displaymath}
From Equation (\ref{7}) of this theorem, we have
\begin{displaymath}
 r(K)=r\left(
  \begin{array}{cccccccccccccc}
             0&   \beta\\
               \beta^{T}&   B\\
  \end{array}
  \right)=r(H)=r(B).
\end{displaymath}
That is
\begin{equation}\label{8}
 r(G)\leq q-1+1+r(H)=q+r(H).
\end{equation}
Combining with  Equations (\ref{6}) and (\ref{8}) of this theorem, we know have $r(G)=q+1+r(H)\leq q+r(H)$, this is a contradiction.
So, $q$ is even.

Let $z$ be a vertex of $C_{q}$ adjacent to $x$. By (a) of Lemma \ref{le:4.1} and Lemmas \ref{le:2.4} and \ref{le:2.5}, we have
\begin{equation}\label{9}
 sr(G^{\sigma})=sr(G^{\sigma}-z)=q-2+sr(K^{\sigma}).
\end{equation}
\begin{equation}\label{10}
 r(G)=r(G-z)+2=q-2+2+r(K)=q+r(K).
\end{equation}

Since $z$ lies on $C_{q}$, so
\begin{equation}\label{11}
 d(G)=d(K)+1=d(H)+1.
\end{equation}

Combining with  Equations (\ref{9})--(\ref{11}) and $G^{\sigma}$ is lower-optimal, we have
$sr(G^{\sigma})=r(G)-2d(G)=q+r(K)-2d(K)-2=q-2+sr(K^{\sigma})$, so
\begin{equation}\label{12}
 sr(K^{\sigma})=r(K)-2d(K).
\end{equation}

 By (a) of Lemma \ref{le:4.1} and Lemmas \ref{le:2.4} and \ref{le:2.5}, we also have
 \begin{equation}\label{13}
 sr(G^{\sigma})=sr(G^{\sigma}-x)=q-2+sr(H^{\sigma}).
\end{equation}
\begin{equation}\label{14}
 r(G)=r(G-x)+2=q-2+2+r(H)=q+r(H).
\end{equation}

Combining with  Equations (\ref{10}) and (\ref{14}), we have
\begin{equation}\label{15}
 r(H)=r(K).
\end{equation}

Combining with Equations (\ref{11}), (\ref{13}) and (\ref{14}), we have
$sr(G^{\sigma})=r(G)-2d(G)=q+r(H)-2d(H)-2=q-2+sr(H^{\sigma})$, so,
\begin{equation}\label{16}
 sr(H^{\sigma})=r(H)-2d(H).
\end{equation}

Combining with Equations (\ref{9}), (\ref{10}), (\ref{12}), (\ref{13}), (\ref{15}) and (\ref{16}), we  obtain (b) and (c) of this theorem.

\textbf{Claim 2}. $q\equiv 2(\rm{mod}~4)$.

Suppose to the contrary that $q=2m$, where $m$ is an even integer. Let $G_{1}=C_{q}-x$, by Lemma \ref{le:2.9}, we have $r(G_{1})=r(C_{q})$. By Lemma \ref{le:2.6}, we have
$r(G)=r(G_{1})+r(G-G_{1})=q-2+r(K)$, which  contradicts to Equation (\ref{10}).

\textbf{Claim 3}. $C_{q}^{\sigma}$ is evenly-oriented.

Suppose to the contrary that $C_{q}^{\sigma}$ is oddly-oriented, by Lemma \ref{le:4.2}, then we have
\begin{equation}\label{17}
 sr(G^{\sigma})=q+sr(H^{\sigma}).
\end{equation}
By Equations (\ref{11}), (\ref{14}) and (\ref{17})and $G^{\sigma}$ is lower-optimal, we have
$$sr(G^{\sigma})=r(G)-2d(G)=q+r(H)-2d(H)-2=q+sr(H^{\sigma}).$$
So,
$sr(H^{\sigma})=r(H)-2d(H)-2$, which contradicts to Equation (\ref{16}).

This completes the proof. \quad $\square$


\noindent\begin{theorem}\label{th:4.4}
Let $y$ be a pendant vertex  of $G^{\sigma}$ adjacent to $x$, and let $H^{\sigma}=G^{\sigma}-y-x$. If $G^{\sigma}$ is lower-optimal, then $x$ does not lie on any cycle of $G$ and $H^{\sigma}$ is lower-optimal.
\end{theorem}
\noindent\textbf{Proof.}
By (c) of Lemma \ref{le:4.1}, we know that $x$ does not lie on any cycle of $G$.

By Lemmas \ref{le:2.4} and \ref{le:2.5}, we have $r(H)=r(G)-2$ and $sr(H^{\sigma})=sr(G^{\sigma})-2$, respectively. Since $x$ does not lie on any cycle of $G$, we have $d(G)=d(H)$. So,

$sr(H^{\sigma})=sr(G^{\sigma})-2=r(G)-2d(G)-2=r(H)+2-2d(H)-2=r(H)-2d(H)$.

This completes the proof. \quad $\square$

\vskip 0.1in

The next paragraph is from \cite{WMT}, which will be useful for later.

\noindent\textbf{In Section 4 of \cite{WMT}}, let $G$ be a graph with pairwise vertex-disjoint cycles, and let $\mathcal{C}(G)$ denote the set of cycles in $G$. By compressing each cycle $O$ of $G$ into a vertex $t_{O}$ we obtain an acyclic graph $T_{G}$ from $G$. More definitely, the vertex set $V(T(G))$ is taken to be $U\cup C_{G}$, where $U$ consists of all vertices of $G$ that do not lie on any cycle and $C_{G}$ consists of vertex $t_{O}$ that is obtained by compressing a cycle $O$, i.e., $C_{G}=\{t_{O}: O\in \mathcal{C}(G)\}$, two vertices in $U$ are adjacent in $T_{G}$ if and only if they are adjacent in $G$, a vertex $u\in U$ is adjacent to a vertex $t_{O}\in C_{G}$ if and only if $u$ is adjacent (in $G$) to a vertex in the cycle $O$, and vertices $t_{O_{1}}$, $t_{O_{2}}$ are adjacent in $T_{G}$ if and only if there exists an edge in $G$ joining a vertex of $O_{1}\in \mathcal{C}(G)$ to a vertex of $O_{2}\in \mathcal{C}(G)$. It is clear that $T_{G}$ is always acyclic. Observe the graph $T_{G}-C_{G}$ (obtained from $T_{G}$ by deleting vertices in $C_{G}$ and the incident edges) is the same as the graph obtained from $G$ by deleting
all cycles and the incident edges, the resultant graph is denoted by $\Gamma_{G}$.

\noindent\begin{theorem}\label{th:4.5}
Let $G^{\sigma}$ be an oriented graph of order $n$. If $G^{\sigma}$ is lower-optimal, then
 \begin{enumerate}[(a)]
  \item  Cycles (if any) of $G^{\sigma}$ are pairwise vertex-disjoint, each cycle $C_{q}^{\sigma}$ of $G^{\sigma}$ is evenly-oriented with order $q\equiv 2(\emph{mod}~4)$.
  \item $r(G)=r(T_{G})+\sum_{O\in\mathcal{C}(G)}|V(O)|$ and $r(T_{G})=r(\Gamma_{G})$.
\end{enumerate}
\end{theorem}
\noindent\textbf{Proof.}
If $G$ has no cycle, then the theorem holds naturally. Suppose $G$ has cycles, let $x$ be a vertex of any cycle. By Lemma \ref{le:4.1}, we know that $x$ lies on only one cycle of $G$, so the first assertion of (a) follows.

We now proceed by induction on the order $n$ to prove the left assertions.

If $n=1$, then all left assertions hold naturally. Suppose the left assertions all hold for any lower-optimal oriented graph of order smaller than $n$, and suppose $G^{\sigma}$ is a lower-optimal oriented graph of order $n\geq2$.

\textbf{Case 1.} If $T_{G}$ has no edges, i.e., $G$ consists of disjoint cycles and some isolated vertices, then the left assertions follow from the following two claims.

\textbf{Claim 1.} $G^{\sigma}$ is lower-optimal if and only if each component of $G^{\sigma}$ is lower-optimal.

\textbf{Claim 2.} A single oriented cycle $C_{q}^{\sigma}$ is lower-optimal if and only if $C_{q}^{\sigma}$ is evenly-oriented with $q\equiv 2(\textrm{mod}~4)$ (by Lemmas \ref{le:2.3} and \ref{le:2.9}).

\textbf{Case 2.} If $T_{G}$ has at least one edge, then $T_{G}$ has at least one pendant vertex $y$. If $y\in U$, then $y$ is also a pendant vertex of $G$. If $y=t_{O}\in C_{G}$, then $G$ has a pendant cycle.

\textbf{Subcase 2.1.} $G$ has a pendant vertex $y$.

Let $x$ be the vertex of $G$ adjacent to $y$, $H^{\sigma}=G^{\sigma}-x-y$. By Theorem \ref{th:4.4}, we know that $x$ is not a vertex of any cycle and $H^{\sigma}$ is lower-optimal. The induction hypothesis to $H^{\sigma}$ implies that
\begin{enumerate}[(1)]
  \item  Each cycle $C_{p}^{\sigma}$ of $H^{\sigma}$ is evenly-oriented with order $p\equiv 2(\rm{mod}~4)$.
  \item $r(H)=r(T_{H})+\sum_{O\in\mathcal{C}(H)}|V(O)|$ and $r(T_{H})=r(\Gamma_{H})$.
\end{enumerate}

Since all cycles of $G$ belong to $H$, we have each cycle $C_{q}^{\sigma}$ of $G^{\sigma}$ is evenly-oriented with order  $q\equiv 2(\rm{mod}~4)$, and $\sum_{O\in\mathcal{C}(H)}|V(O)|=\sum_{O\in\mathcal{C}(G)}|V(O)|$. Noting that $y$ is also a pendant vertex of $T_{G}$ (resp., $\Gamma_{G}$) adjacent to $x$ and $T_{H}=T_{G}-x-y$ (resp., $\Gamma_{H}=\Gamma_{G}-x-y$), combining with (2) of Subcase 2.1 and Lemma \ref{le:2.4}, then we have
$$r(G)=r(H)+2=r(T_{H})+\sum_{O\in\mathcal{C}(H)}|V(O)|+2=r(T_{G})+\sum_{O\in\mathcal{C}(G)}|V(O)|,$$
and
$$r(T_{G})=r(T_{H})+2=r(\Gamma_{H})+2=r(\Gamma_{G}).$$

\textbf{Subcase 2.2.} $G$ has a pendant cycle $C_{q}$.

Let $x$ be the unique vertex of $C_{q}$ of degree 3, $H^{\sigma}=G^{\sigma}-C_{q}^{\sigma}$ and $K^{\sigma}=H^{\sigma}+x$. By (c) of Theorem \ref{th:4.3}, we know that both $H^{\sigma}$ and $K^{\sigma}$ are lower-optimal. The induction hypothesis to $K^{\sigma}$ implies that
\begin{enumerate}[(i)]
  \item  Each cycle $C_{p}^{\sigma}$ of $K^{\sigma}$ is evenly-oriented with order $p\equiv 2(\rm{mod}~4)$.
  \item $r(K)=r(T_{K})+\sum_{O\in\mathcal{C}(K)}|V(O)|$ and $r(T_{K})=r(\Gamma_{K})$.
\end{enumerate}

Combining with (a) of Theorem \ref{th:4.3}, assertion (i) of Subcase 2.2 and $\mathcal{C}(G)=\mathcal{C}(K)\cup{\{C_{q}\}}$ imply that each cycle of $G^{\sigma}$ is evenly-oriented with order  $q\equiv 2(\rm{mod}~4)$. Applying (b) of Theorem \ref{th:4.3} and  assertion (ii) of Subcase 2.2, we have
\begin{equation}\label{18}
 r(G)=q+r(K)=q+r(T_{K})+\sum_{O\in\mathcal{C}(K)}|V(O)|.
\end{equation}
Since $T_{K}$ is isomorphic to $T_{G}$, and $q+\sum_{O\in\mathcal{C}(K)}|V(O)|=\sum_{O\in\mathcal{C}(G)}|V(O)|$, it follows from Equation (\ref{18}) that
\begin{equation}\label{19}
 r(G)=r(T_{G})+\sum_{O\in\mathcal{C}(G)}|V(O)|,
\end{equation}
which proves the first assertion of (b) of this theorem.

By (b) of  Theorem \ref{th:4.3}, we have
\begin{equation}\label{20}
 r(G)=q+r(H).
\end{equation}
Noting that $\mathcal{C}(G)=\mathcal{C}(H)\cup{\{C_{q}\}}$, then from Equations (\ref{19}) and (\ref{20}), we have
\begin{equation}\label{21}
 r(T_{G})=r(G)-\sum_{O\in\mathcal{C}(G)}|V(O)|=q+r(H)-\sum_{O\in\mathcal{C}(G)}|V(O)|=r(H)-\sum_{O\in\mathcal{C}(H)}|V(O)|.
\end{equation}
Since $H^{\sigma}$ is also lower-optimal, the first assertion of (b) of this theorem applying to $H$ implies that
\begin{equation}\label{22}
 r(H)=r(T_{H})+\sum_{O\in\mathcal{C}(H)}|V(O)|.
\end{equation}
Equations  (\ref{21}) and (\ref{22}) implies that
\begin{equation}\label{23}
 r(T_{G})=r(T_{H}).
\end{equation}
The induction hypothesis to $H^{\sigma}$ implies that
\begin{equation}\label{24}
 r(T_{H})=r(\Gamma_{H}).
\end{equation}
Since $\Gamma_{G}=\Gamma_{H}$, combining with Equations  (\ref{23}) and (\ref{24}), we have $r(T_{G})=r(\Gamma_{G})$.

This completes the proof. \quad $\square$

Let $T$ be an acyclic graph with at least one edge, we denote by $\widetilde{T}$ the subgraph obtained from $T$ by deleting all pendant vertices of $T$.
\noindent\begin{lemma}\label{le:4.6}(\cite{MWTLAA})
 Let $T$ be an acyclic graph with at least one edge. Then
 \begin{enumerate}[(a)]
  \item  $r(\widetilde{T})<r(T)$.
  \item If $r(T-W)=r(T)$ for a subset $W$ of $V(T)$, then there is a pendant vertex $v$ such that $v\notin W$.
\end{enumerate}
\end{lemma}

\noindent\textbf{Proof of Theorem \ref{th:1.3}.}

\textbf{Sufficiency:}
Suppose that $G^{\sigma}$ meets all the conditions (1)--(3) in Theorem \ref{th:1.3} and $k$ steps of $\delta$-transformations can switch $G$ to a crucial subgraph of $G_{0}$, which is the disjoint union of $d(G)$ cycles together with $l$ isolated vertices. By Lemmas \ref{le:2.4} and \ref{le:2.5}, we have
\begin{equation}\label{25}
 sr(G^{\sigma})=2k+sr(G_{0}^{\sigma}),~r(G)=2k+r(G_{0}).
\end{equation}

Since each cycle $C_{q}$ of the crucial subgraph $G_{0}$ of $G$ is evenly-oriented with order $q\equiv 2(\rm{mod}~4)$, by Lemmas \ref{le:2.3} and \ref{le:2.9}, we have
\begin{equation}\label{26}
 sr(G_{0}^{\sigma})=\sum_{O\in\mathcal{C}(G)}sr(O^{\sigma})=\sum_{O\in\mathcal{C}(G)}|V(O)|-2d(G)=\sum_{O\in\mathcal{C}(G)}r(O)-2d(G)=r(G_{0})-2d(G).
\end{equation}
By Equalities (\ref{25}) and (\ref{26}), we have
$$sr(G^{\sigma})=2k+sr(G_{0}^{\sigma})=2k+r(G_{0})-2d(G)=2k+r(G)-2k-2d(G)=r(G)-2d(G).$$

This completes the proof of sufficiency.

\textbf{Necessity:} Let $G^{\sigma}$ be a lower-optimal oriented graph. By (a) of Theorem \ref{th:4.5}, we can obtain the (1) and (2) of Theorem \ref{th:1.3}. Thus $G$ has precisely $d(G)$ vertex-disjoint cycles, and the acyclic graph $T_{G}$ respect to $G$ is well defined. Now, we will proceed by induction on the order $n$ of $G^{\sigma}$ to prove (3) of Theorem \ref{th:1.3}.

\textbf{Case 1.} If $n=1$, then the assertion holds naturally.

\textbf{Case 2.} Suppose the assertion holds for all lower-optimal oriented graphs with order smaller than $n$, and let $G^{\sigma}$ be a lower-optimal oriented graph of order $n$.

\textbf{Subcase 2.1.} If $T_{G}$ has no edges, then $G$ is the disjoint union of $d(G)$ cycles along with some isolated vertices, and the assertion holds naturally.

\textbf{Subcase 2.2.} If  $T_{G}$ has at least one edge, by (b) of Theorem \ref{th:4.5}, we have $$r(T_{G})=r(\Gamma_{G})=r(T_{G}-C_{G}).$$ (b) of Lemma \ref{le:4.6} shows that there is a pendant vertex of $T_{G}$ not in $C_{G}$. Thus $G$ has at least one pendant vertex. Let $y$ be a pendant vertex of $G$ adjacent to a vertex $x$ of $G$, by Theorem \ref{th:4.4}, $x$ does not lie on any cycle of $G$ and the graph  $H^{\sigma}=G^{\sigma}-x-y$  is also lower-optimal, and also has $d(G)$ cycles. The induction hypothesis applying to $H^{\sigma}$ implies that a series of $\delta$-transformations can switch $H$ to a crucial subgraph of $G_{0}$ consisting of $d(G)$ disjoint union  cycles together with some isolated vertices. Combining with the first step of $\delta$-transformation applying to $G$ and all the other $\delta$-transformations done latter, we can switch $G$ to the crucial subgraph $G_{0}$.

This completes the proof. \quad $\square$


\begin{thebibliography}{99}



\bibitem{BBD} J.H. Bevis, K.K. Blount, G.J. Davis, The rank of graph after vertex addition. Linear Algebra Appl. 265(1997) 55--69.

\bibitem{BM}
 J.A. Bondy,  U.S.R. Murty, Graph Theory with Applications. Elsevier. New York (1976)

\bibitem{CT} Li Chen, Fenglei Tian, Skew-rank of an oriented graph with edge-disjoint cycles. Linear and Multilinear Algebra. 64(2016) 1197--1206.

\bibitem{CL} Bo Cheng, Bolian Liu, On the nullity of graphs.  Electron. J. Linear Algebra 16(2007) 60--67.

\bibitem{CS} L. Collatz, U. Sinogowitz, Spektren endlicher grafen. Abh. Math. Sem. Univ. Hamburg. 21(1957) 63--77.


\bibitem{GX} Shicai Gong, Guanghui Xu, On the nullity of a graph with cut-point. Linear Algebra Appl. 436(2012) 135--142.


\bibitem{LXYG} Xueliang Li, Guihai Yu,  The skew-rank of oriented graphs. Sci. Sin. Math. 45(2015) 93--104. (in Chinese)

\bibitem{LWZ} Yong Lu, Ligong Wang, Qiannan Zhou, Bicyclic oriented graphs with skew-rank 6. Appl. Math. Comput.  270(2015) 899--908.




\bibitem{MWTDAM} Xiaobin Ma,  Dein Wong, Fenglei Tian, Nullity of a graph in terms of the dimension of cycle space and the number of pendant vertices. Discrete Appl. Math. 215(2016) 171--176.


\bibitem{MWTLAA} Xiaobin Ma,  Dein Wong, Fenglei Tian, Skew-rank of an oriented graph in terms of matching number. Linear Algebra Appl. 495(2016) 242--255.


\bibitem{QHYG} Hui Qu, Guihai Yu, Bicyclic oriented graphs with skew-rank 2 or 4. Appl. Math. Comput.  258(2015) 182--191.

\bibitem{QHYGFL}Hui Qu, Guihai Yu, Lihua Feng, More on the minimum skew-rank of graphs. Oper. Matrices. 9(2015) 311--324.

\bibitem{RCZ} S. Rula, An Chang, Yirong Zheng, The extremal graphs with respect to their nullity. J. Inequal. Appl. 2016(2016) 65


\bibitem{SST} Yazhi Song, Xiaoqiu Song, Bit-Shun Tam, A characterization of graphs $G$ with nullity $V(G)-2m(G)+2c(G)$. Linear Algebra Appl. 465(2015) 363--375.

\bibitem{LWANG} Long Wang, Characterization of graphs with given order, given size and given matching number that minimize nullity. Discrete Math. 339(2016) 1574--1582.


\bibitem{WW} Long Wang, Dein Wong, Bounds for the matching number, the edge charomatic numbber and the independence number of a graph in terms of rank. Discrete Appl. Math. 166(2014) 276--281.

\bibitem{WMT}
Dein Wong, Xiaobin Ma, Fenglei Tian, Relation between the skew-rank of an oriented graph and the rank of its underlying graph. European J. Combin. 54(2016) 76--86.




\end{thebibliography}
\end{document}